\documentclass[12pt]{article}
\usepackage{MnSymbol}
\usepackage{amsmath,amsfonts,
verbatim}
\usepackage[all,cmtip]{xy}

\newcommand{\Q}{{\mathbb Q}}
\newcommand{\C}{{\mathbb C}}
\newcommand{\N}{{\mathbb N}}
\newcommand{\Z}{{\mathbb Z}}
\newcommand{\R}{{\mathbb R}}

\newcommand{\F}{\mathrm{F}}
\newcommand{\U}{{\mathbb U}}

\newcommand{\D}{\mathbb{D}}

\newcommand{\pp}{\mathfrak{p}}
\newcommand{\lfr}{\mathfrak{l}}
\newcommand{\ii}{\mathfrak{i}}

\newcommand{\ZZ}{{^*\mathbb{Z}}} 
\newcommand{\la}{\langle}
\newcommand{\ra}{\rangle}
\newtheorem{pkt}{}[section]  
\newcommand{\bpk}{\begin{pkt}\rm }  
\newcommand{\epk}{\end{pkt}} 
\newcommand{\inv}{^{-1}}   
\newcommand{\dd}{\mathbf{d}}
\newcommand{\be}{\begin{equation}}  
\newcommand{\ee}{\end{equation}}  
\newcommand{\ml}{\mathsf{lm}}

\newcommand{\K}{\mathrm{K}}

\newcommand{\uu}{\mathbf{u}}
\newcommand{\vv}{\mathbf{v}}
\newcommand{\II}{\mathrm{I}}

\usepackage{graphicx}

\usepackage{tikz}

\title{Physics over a finite field and Wick rotation}
\author{B.Zilber}

\begin{document}
\maketitle
\begin{abstract}
The paper develops some mathematics supporting an earlier hypothesis that the physical universe is a finite system co-ordinatised by a huge finite field $\F_\pp$ which looks like the field of complex numbers to an observer.

Earlier we constructed a place ('limit' homomorphism) $\ml$ from a pseudo-finite field $\F_\pp$ onto the compactified field of complex numbers. In the current paper we construct  $\ml$ in a more concrete form. In particular, $\ml$ sends certain multiplicative subgroups $'\R'_+$ and $'\mathbb{S}'$ of $\F_\pp$  onto the non-negative reals $\R_+$ and the unit circle  $\mathbb{S}$ in $\C.$    Thus $\F_\pp,$ $'\R'_+$ and $'\mathbb{S}'$ provide co-ordinates for physical universe.

We introduce two systems of natural units corresponding to $'\R'_+$ and $'\mathbb{S}',$ respectively, on the logarithmic scale.
 The passage from the scale of units of $'\R'_+$ to the scale of units of $'\mathbb{S}'$ corresponds to a multiplication (on the logarithmic scale) by a 'huge' (non-standard) integer $\ii$ equal approximately  to $\sqrt{\pp}.$ This provides an explanation to the phenomenon of  Wick rotation.

In the same model we explain the phenomenon of phase transition in a large finite system     

\end{abstract}
\section{Introduction}
\bpk
The hypothesis that the universe is infinite is an open question. This concerns both the size of the universe and the number of atoms or elements that comprise it. Since it is now accepted that there is a minimal length, the Planck length, assumption of the spatial finiteness of the universe implies the assumption of the finiteness of the number of its elements. 

In \cite{Perfect} we discussed the concept of approximation in physics and  the suggestion that  physics universe is co-ordinatised by a {\em huge}\footnote{We use {\em huge} for numbers which some physics authors call ``ridiculously large'', see  the discussion in \ref{Ds3} below.}
 finite field $\F_\pp.$ It was proved (Proposition 5.2 of \cite{Perfect}) that the only metric field (locally compact field) that can be approximated by finite fields is the field $\C$ of complex numbers. Thus ``seen from afar'' the huge finite field looks like a field of complex numbers $\C.$

The current work has been inspired by Hao Hu, the expert in Philosphy of Physics, who approached the author with the suggestion to apply  the idea of \cite{Perfect} to statistical mechanics and to
attack the well-discussed problem of phase transition (which happen in large finite systems but require the assumption of infinity for its mathematical theory). We suggest here an answer to the problem. 

Perhaps a more important outcome of the mathematical theory developed below is the interpretation of the phenomenon of  {\em Wick rotation} as the result of the change of scales in physics.

\medskip

I would also like to note that a compatible  attempt to develop the mathematical background for a theory of finite physics was presented in  
 \cite{QM0}, \cite{QM1} and \cite{PaperN}.
\epk
\bpk Let us recall the notion of structural approximation suggested in \cite{Perfect}, in the specific context of approximation by huge finite fields $\F_\pp.$  A finite structure is discrete but to see its grainy structure we should be able to detect difference, {\em inequality}, between its neighbouring elements, which might be impossible with the instruments  we use to observe the structure. However, if there is a shape within the $n$-space $\F_\pp^n$ which is given by an algebraic {\em equation} then the observer could see it as a shape in a continuous field, say equal to $\R,$ $\C$ or  maybe a p-adic field, given by the same equation. This brings us to the idea that such an approximation by an observer is a map (called ``limit'') onto a continuous field $\K,$
$$\ml: \F_\pp\to \K$$
which takes tuples $(x_1,\ldots,x_n)\in \F_\pp^n$ satisfying a polynomial equation\linebreak $f(x_1,\ldots,x_n)=0$ (with integer coefficients)
to the tuple $(y_1,\ldots,y_n)\in \K^n$ satisfying the same equation. In other words, $\ml$ is a ring-homomorphism. 

In fact, for a finite $\F_\pp$ this scheme is not going to work verbatim  but it works when we assume that $\F_\pp$ is infinite {\em pseudo-finite}, which for all intents and purposes replaces  a huge finite structure, see \ref{1.3} below for  definition. 

However, as explained in \cite{Perfect}, the requirement that $\ml$ is defined on the whole of the discrete structure necessitates that the target structure $\K$ must be compact, which for a metric field can be achieved by adding a point $\infty$ (so for $\K=\C$ the compactification gives us the {\em extended complex numbers} $\bar{\C}:=\C\cup \{ \infty\},$ equivalently, the  Riemann sphere ). In particular, there are non-zero elements $x\in \F_\pp$ indistinguishable from $0$ (that is $\ml(x)=0$) and so for the inverses $x\inv \in \F_\pp,$  $\ml(x\inv)=\infty.$  Such a map $\ml$ between fields is called a {\em place} in algebra. 

\medskip

The above mentioned key Proposition 5.2 of \cite{Perfect} states: 

{\em There exists a place \be \label{place1}\ml: \F_\pp\to \bar{\C} \ee
 from any pseudo-finite field $\F_\pp$ of characteristic $0$ onto the compactification of the field of complex numbers $\C,$ and $\C$ is the only metric (locally compact) field for which such a map exists. }   

\epk
\bpk {\bf The relationship with the p-adic approach in physics}. The p-adic approach in physics, in particular string theory, has proved quite productive, see e.g. the survey \cite{Chren}. 
The recourse to a prime $p$ is motivated by the needs of discretisation and the field $\Q_p$ of $p$-adic numbers has the advantage of bearing a nice metric and locally compact topology. There is no preferred prime but a very large prime like the above $\pp$ seems to be a reasonable choice. We note that by definition there is a canonical place
\be \label{place2}\Q_\pp\to \F_\pp\ee and thus, combining with (\ref{place1}) we get a place
$$\Q_\pp\to \bar{\C}.$$
In other words,  p-adic calculations pass via (\ref{place2}) to $\F_\pp$ which, according to (\ref{place1}), can be reinterpreted as calculations in  the complex numbers.
\epk
\bpk The main mathematical
 problem in making a practical use of the idea of a ``physics over a  finite field $\F_\pp$''  is to find a way of representing the common-sense real quantities of physics inside the huge finite field $\F_\pp$ 
 and to find such a representation that allows the standard physics calculations. That is to explain how both $\R$ and $i\R$ emerge from a large finite field.
 
In this regard it is useful to invoke a notion of {\em feasible numbers} that was discussed by philosphers of mathematics, mathematical logicians and computer scientists, see e.g. \cite{Sazonov} for a mathematical treatment of the notion. Roughly speaking, $1,2,\ldots,1000,$ as well as
their ratios such as $1/5$ and  $0.203$ are feasible numbers, but   the Avogadro number $\sim 10^{23}$ is not feasible. We think of the latter as {\em very large} but potentially observable numbers. This contrasts with {\em huge } numbers, such as $\pp,$ the number of points in  $\F_\pp,$ and $\pp>> 10^{23}.$   

Thus, while we work with feasible numbers $1,2,\ldots $  inside $\F_\pp$  we can think of these as the usual  integers 
 but when our integers become  very large but still much less than $\pp$ then  $\ml$ takes such ones to $\infty.$ Further on this scale,  according to the approximation theorem of \cite{Perfect}, the integers start to behave like complex numbers, e.g. if an integer $\ii$ satisfies $\ii^2+1=0 \mod \pp$ (huge number)
it takes the role of $\sqrt{-1}.$ 

All this is made precise in the formulation of the Main Theorem in  \ref{MainA} below. 
\epk
\bpk\label{1.3} In the current work 
a ``huge finite field'' is a {\em pseudo-finite field} $\F_\pp$ which can be obtained by considering a non-principal ultrafilter $\mathcal{D}$ on the set of prime numbers $P\subset \N$ (positive integers) and the ultraproduct
$$\F_\pp=\prod_{p\in P}\F_p/\mathcal{D}.$$
Such a construction sees $\F_\pp$ as a logical limit of finite fields $\F_p$ along the ultrafilter: a first order sentence $\Phi$ is valid in $\F_\pp$ if and only if it is valid in almost all, in the sense of $\mathcal{D},$ finite fields $\F_p.$

Model theory tells us that $\F_\pp$ can equaivalently be  obtained as the quotient-ring of the ring $\ZZ$ of non-standard integers by the prime ideal $\pp\, \ZZ,$ where $\pp$ is the respective non-standard prime number,
$$ \F_\pp\cong \ZZ/\pp\mbox{ where } \ZZ=\Z^P/\mathcal{D}.$$

The interpretation of $\F_\pp$ in non-standard integers  allows us to apply, among others, the means of non-standard analysis, in particular the {\em standard part map}
$$\mathrm{st}: {^*\Q}\to \R\cup \{ +\infty,\, -\infty\}$$
from non-standard rationals $\frac{l}{m},$ $l,m\in \ZZ,$ to the compactification of the reals. 
   
\epk

\bpk\label{1.4} Along with $\pp$ and the field $(\F_\pp;+,\cdot,0, 1)$ we specify:

-  a non-standard {\em higly divisible} number $\lfr$ (each standard integer  $m$ divides $\lfr$) satisfying  some other assumptions below;

- a two-sorted pseudo-finite structure $(\U_{\pp,\lfr},\F_\pp)$ with 
$$\U_{\pp,\lfr}=(\ZZ/_{(\pp-1)\mathfrak{l}};+,\hat{0},\hat{1})$$ a pseudo-finite additive cyclic\footnote{here and below ``cyclic'' means in the pseudo-finite sense, i.e. the ultraproduct of cyclic groups } group  of order $(\pp-1)\mathfrak{l}$ with generator $\hat{1};$

- a surjective   group homomorphism
$$\exp_\pp: \U_{\pp,\lfr}\to \F_\pp^\times;  \ \  \ n\cdot \hat{1}\mapsto \epsilon^n$$
where $\epsilon$ is a generator of the (pseudo)-cyclic group $\F_\pp^\times,$ $n\in \ZZ.$
It follows that $$\ker \exp_\pp=(\pp-1)\cdot \U_{\pp,\lfr},$$
the subgroup generated by $(\pp-1)\cdot\hat{1},$ (and suggests that $\pp-1$ of $\U_{\pp,\lfr}$ should be interpreted as $2\pi i$);



- a pair of surjective ``limit''  homomorphisms (places) $\ml$  which make the diagram commute
\be\label{Diag1} \begin{array}{lllll}
\ml_\U:\ \ \ \U_{\pp,\lfr}\ \ \ \ \twoheadrightarrow \ \ \ \ \ 
\bar{\C}\\
\\

\ \ \ \exp_\pp \downarrow \ \ \ \ \ \ \ \ \ \ \ \ \exp \downarrow\\

\\

\ml_\F: \ \ \ \F_\pp\ \ \ \ \to \ \ \ \ \ \ \ \ \bar{\C}  
\end{array}    \ee

There is a natural cyclic order on $\U_{\pp,\lfr}$ in which $u+\hat{1}>u$ and there is a related cyclic order on $\F_\pp^\times$ in which $\epsilon^{n+1}>\epsilon^n$ for all $n\in \ZZ.$

We treat $\hat{1}$ as an infinitesimal and choose two ``units of length'' $\uu$ and $\vv,$ elements of $\U_{\pp,\lfr},$ 
$$\hat{1}<< \uu<<\vv.$$
More precisely,
\be\label{defuv} \uu=\frac{\pp-1}{\ii}\mbox{ and }\vv=\pp-1\ee
for some $\ii\in {^*\N}$ such that  \be\label{il} \ii\lfr|(\pp-1),\ \ \ii>>\lfr\ee
($\ii \lfr$ divides $\pp-1$).

Thus, $\uu$ and $\vv$ are units of two very different scales (see \ref{concl1} below for further comment).

 We  assume that \be\label{lambdamu} \ii=\iota^2\mbox{ and }\lfr=\mu^2\ee for some $\mu,\iota\in {^*\N}.$ 

We also need to assume \be \label{ii} \ii^2+1=\pp \mbox{ \ \ or\ \ \  $\ii,\lfr$  algebraically independent in }\F_\pp\ee

(the first is the preferable and more elegant assumption but it is not known 
whether it is consistent with $\pp$ being infinite (the Landau problem)). 

It is easy to check that our set of assumptions (\ref{il})-(\ref{ii}) along with the assumption of high divisibility of $\lfr$ are consistent. 
We are going to slightly extend these assumption later, in particular (\ref{ccmax}) assumes that $\lfr^n<\pp$ for all $n\in \N.$  

\medskip

We also define additive subgroups  of $\U_{\pp,\lfr}$ called suggestively
$'\R'$ and $'i\R'.$

\epk
\bpk\label{MainA} {\bf Main Theorem.}  {\em
There exists a  surjective ring homomorphism (place)
$$\ml_\F: \F_\pp\twoheadrightarrow \bar{\C}$$
and a surjective additive semigroup homomorphism 
$$\ml_\U: \U_{\pp,\lfr} \twoheadrightarrow\bar{\C}$$
such that:  

 the diagram (\ref{Diag1}) is commutative,  
\be \label{R} \ml_\U: {'\R'}\twoheadrightarrow \R; 
\ee
\be \label{iR} \ml_\U: {'i\R'}\twoheadrightarrow i\R; 
\ee
where $i=\sqrt{-1}.$

For  any $l\in \ZZ$ such that $0<l\le \lfr,$ 

\be\label{ll} \ml_\F: l\cdot \lfr\inv \mapsto  \mathrm{st}(\frac{l}{\lfr})\in \R
\ee

\be \label{e2} \ml_\F: 
\exp_\pp(l\uu)\mapsto e^{-\pi z},\ \  \mbox{ for } z:=\mathrm{st}
(\frac{l}{\lfr})\in  \R
\ee
\be \label{e1} \ml_\F: 
\exp_\pp(l\vv)\mapsto e^{i\pi z}, \ \ \mbox{ for } z:=\mathrm{st}
(\frac{l}{\lfr})\in  \R\ee
where $\mathrm{st}$ is the standard part map.

For $a\in \Q_+,$ a complete square and $\mu\in \ZZ$ such that $\mu^2=\lfr:$ 
\be \label{integral2} \ml_\F:\
 \frac{1}{\mu}\sum_{-\ii\lfr/a\le n<\ii\lfr/a} \exp_\pp( a\frac{n^2}{2\lfr}\uu)\ \ \mapsto \ \
 \int_\R e^{-a\pi x^2} dx
\ee
and $\int_\R e^{- a x^2} dx= \frac{1}{\sqrt{a}},$ according to the standard definition  of the Gaussian integral;

 \be \label{integral1} \ml_\F:\  \frac{1}{\mu}\sum_{-\lfr/a \le n<\lfr/a} \exp_\pp( a\frac{n^2}{2\lfr}\vv)\ \ \mapsto\ \ \int_\R e^{ia\pi  x^2} dx \ee
where $\int_\R e^{ia\pi   x^2} dx= \frac{e^\frac{\pi i}{4}}{\sqrt{a}}$ according to  the Quantum Mechanics calculus.}
 \medskip
 
Note that (\ref{R}) and (\ref{iR})  give us subgroups of $\F_\pp^\times$
$${'\R'}_+:=\exp_\pp('\R')\mbox{ and }{'\mathbb{S}'}:=\exp_\pp('i\R')$$
 furnish a good analogue of polar coordinate system in $\F_\pp.$
\epk
\bpk {\bf Remark.} Since $\ml$ is a place, (\ref{ll}) determines the values of $\ml(Q(l)$ for any rational function $Q(x)$ over $\Z.$  
\epk
\bpk\label{concl1} {\bf Discussion}. The theorem clarifies the relationship between  two different scales in $(\U_{\pp,\lfr},\F_\pp)$ presented by units $\uu$ and $\vv.$ 
These should be though of as units for physics of 'low energy' and 'high energy', respectively, the latter being physics at quantum level and the former the physics at the level of  Brownian motion.
 
The statements (\ref{e2}) and (\ref{e1}) demonstrate that the action of the huge pseudofinite integer $\ii$ which changes the scale of units in $\U_{\pp,\lfr}$ (recall that $\vv=\ii\cdot \uu$, and  $\ii\approx \sqrt{\pp}$) is seen in $\C$ as the {\bf  Wick rotation} $e^{\pi z}\mapsto e^{i\pi z}.$

 The Gaussian integrals in (\ref{integral2}) and (\ref{integral1}) are mathematical manifestation of the same phenomena expressed by the summation formulae over $\F_\pp$ on the left hand side of $\mapsto.$ The difference in the integral expressions on the right  comes from the difference in the scale of units  that measure the respective processes.

Thus the puzzling regularity of the transition from the Brownian motion integral (\ref{integral2}) to the quantum mechanics integral (\ref{integral1}) known and exploited by physicists as   Wick rotation has an explanation as a mathematical consequence of the change of scales. 
\epk
\bpk {\bf Discussion.} The ``limit''  $\ml_\F$ maps the discrete field $\F_\pp$ into (the compactification of) the field of complex numbers and thus endows the image with a metric. The Theorem determines $\ml_\F$
on the subring generated by  specific points (see (\ref{ll}), (\ref{e2}) and (\ref{e1})) but leaves the rest free. This means that the observer, which sees $\F_\pp$ with all the algebraic geometry over it through $\ml_\F,$ has some freedom in choosing the metric on algebraic varieties. In particular, if $Z\subseteq \F_\pp^m$ is an algebraic subvariety, say a torus, then
$\ml_\F(Z)$ is a subset of the compactification of $\C^m,$ a compact complex variety.  In general, such a compactification is far of being unique.

Thus the freedom in the choice of $\ml_\F$ implies a respective degree of freedom in the choice of complex/metric version of physics. 
\epk

\section{ Statistical physics and phase transition}
\bpk {\bf Physical units and dimensions}

In the formalism of two sorts $\U_{\pp,\lfr}$ and $\F_\pp$ introduced in \ref{1.4} the sort $\U_{\pp,\lfr}$ is assumed to be the sort thar holds all the physical units (dimensions). It is convenient for each principal unit of measurement 
to define a special sort $\D_i,$ $i=1,\ldots,k$ which is going to be naturally interpreted in terms of $\U_{\pp,\lfr}.$
.

Each sort $\D_i$ is a subgroup of $\U_{\pp,\lfr}$ and so has a cyclic additive group structure with the unit (generator) $\dd_i\in \U_{\pp,\lfr}$ ($\D_i$-unit). $\D_i$ is isomorphic to  $\U_{\pp,\lfr}/\ker_i,$ where 
$$\ker_i=\frac{(\pp-1)\lfr}{\dd_i}\cdot \U_{\pp,\lfr},\ \ \dd_i|\,(\pp-1)\lfr.$$
Thus the size of $\D_i,$
$$|\D_i|=\frac{(\pp-1)\lfr}{\dd_i}.$$
  
 Between some of the unit sorts there are bilinear  maps
$$\D_1\times \D_2\twoheadrightarrow \D_{3}; \ \ \ (x_1\dd_1,x_2\dd_2)\mapsto x_1x_2\dd_3$$
where we assume $\ker_3= \ker_1\cap\ker_2$ and
 $$x_1=u_1+\ker_1,\ x_2=u_1+\ker_2,\ \
\ x_1\cdot x_2:= u_1\cdot u_2+\ker_3$$  in the ring $\ZZ/(\pp-1)\lfr.$   

Thus counting in units of $\D_3$ gives
$$x_3=x_1x_2$$
and this can be equivalently written as
$$x_1=x_2\inv x_3.$$

The units $\uu$ and $\vv$ introduced in \ref{1.4} are examples of $\D$-units.
\epk
\bpk The field sort $\F_\pp$ is assumed to be dimensionless and the exponentiation map $\exp_\pp$ restricted to a sort $\D_i$ is a homomorphism
$$\exp_\pp: \D_i\twoheadrightarrow \F_\pp^\times; \ \ n\cdot \dd_i\mapsto \exp_\pp(n\dd_i).$$
\epk
\bpk In Statistical Mechanics the dimensions in $\U_{\pp,\lfr}$ are usually  Energy (E), Temperature (T) and (in ferromagnets) magnetic moment (H).

According to this theory probability that the system in temperature $T$ is in a state $\sigma$  is equal to $\exp(-\frac{E_\sigma}{kT}),$ and the probability $p_n$ that the system consists of exactly $n$ atoms, out of possible $N,$ is
$$p_n=\frac{P_n}{Z_N},\ \ P_n=\sum_{\sigma\in \Sigma(n)} \exp(-\frac{E_\sigma}{kT}),$$
where $\Sigma(n)$ are all the states with exactly $n$ atoms and
$$Z_N=  \sum_{\sigma} \exp(-\frac{E_\sigma}{kT}),$$
where $\sigma$ runs in all possible states with at most $N$ atoms.

Setting (with some simplifications) $y:=\exp(-\frac{H}{kT}),$
the equilibrium state of the system of volume $N$ (that is having up to $N$ particles) is analysed via the polynomial 
 $$\mathcal{P}_N(y):=\sum_{n=0}^N p_n y^n.$$
 which is called  the {\em grand partition function} of the system.
 
 Assuming frequentist probabilities  
we may regard the $P_n$ and $Z_N$  integers.
  
Since there are very few  restrictions on states in the models the number $Z_N$ is close to the number of all possible subsets, that is
\be\label{Z} \sum _{n=0}^N P_n=Z_N\approx 2^N.\ee

(The estimate (\ref{Z}) appears also in \cite{Tong} on page 4, for $N$ the Avogadro number, and the number $Z_N$ is being characterised as {\em ridiculously large number}.)

Also note that when $n$ is near $N/2,$ $P_n$ is near its maximum 
 \be \label{Pn} P_n\approx \left(\begin{array}{ll}  N\\ \frac{N}{2}
 \end{array}\right)\approx \frac{2^N}{\sqrt{N}}\ee

\epk
\bpk 
We note that since $N$ is supposed to be very large number, by (\ref{Pn}) $P_n$ can reach huge values and so  modulo $\pp$ will be outside the real part $'\R'$ of $\F_\pp.$ Same is true for $Z_N.$ Hence in general the $p_n$ should be treated as {\bf probability amplitudes} rather than classical probabilities, and $\mathcal{P}_N$ should be treated as a polynomial with complex coefficients rather than polynomial over $\R.$     

\epk
\bpk The seminal work of C.N.Yang and T.D.Lee,  \cite{YL1} - 
\cite{YL2} (1952), laid the ground for the modern theory of critical points in the evolution of large finite systems such as the ideal gas. 

The main theorem of \cite{YL1} states that the phase transition in the system happens at the point $y_\mathrm{crit}$ (a critical point) such that
$$   \mathcal{P}_N(y_\mathrm{crit})=0$$
which obviosly {\bf can not be a real point}. 
The paper
   analyses complex roots of the polynomial which proves very important for the behaviour of the system near the critical point.

The modern theory resolves this paradox by assuming $N\to\infty,$ in which case  $\frac{1}{N}\ln\mathcal{P}_N(y)$ converges, away of the critical point, to an analytic function (the passage to the thermodynamical limit)
and $y_\mathrm{crit}$ converges to a real point. This solution of the paradox is not considered to be fully satisfactory as the actual systems are always finite, although very large. Under the passage to thermodynamical limit some information is being lost.
\epk
\bpk

 The hypothesis of the {\em universe over finite field} suggests a solution to the paradox. Under the hypothesis $y_\mathrm{crit}\in \F_\pp,$  that is $y_\mathrm{crit}$ is an integer such that

\be\label{V} \mathcal{P}_N(y_\mathrm{crit})=0\mod \pp.\ee 

One can explore the assumption further using the reduction (24) from \cite{YL2} together with Theorem 3 therein which states that after the reduction all the zeroes of  $\mathcal{P}_N$ are on the unit circle and the limit of the zeroes is $=1.$ This leads us to the conclusion that
 $$\mathcal{P}_N(1)=0\mod \pp$$

 Equivalently, $$\sum_{n=1}^N P_n=0\mod \pp.$$
 
 Combining this with (\ref{Z}) and assuming that
   the number of atoms $N$ is $\approx 10^{23},$ the Avogadro number,
we can make an estimate on $\pp,$ the upper bound:  \be\label{Avogadro} \pp < 2^N\approx 2^{10^{23}}\ee

On the other hand the same argument proves that {\em there is a low bound on the volume $N$ of gas which  allows a phase transition}, that is has a critical point  satisfying (\ref{V}):
$$N> \log \pp.$$
 
\epk
\bpk\label{Ds3} {\bf Discussion.} The assumption of physics over a finite field explains the necessity of extending the definition of grand partition function as the function of complex variable as well as explains  the phase transition in a large finite system. 

\medskip

In a forthcoming paper we  develop an {\em analytic theory } on $\U_{\pp,\lfr}$ which, via $\ml_\U$ corresponds to the analytic theory on $\C.$ In particular the expression like 
$$\frac{1}{N}\ln\mathcal{P}_N(y) \mbox{ and }\ml_\U \left\lbrace \frac{1}{N}\ln\mathcal{P}_N(y)\right\rbrace$$
become the lawfull objects of the theory and one can carry on the thermodynamic theory as usual, along with the discrete theory on $\F_\pp.$ 
 
\epk

The rest of the paper is purely mathematical. It provides the  construction of the limit map $\ml$ and proofs.

\section{The pseudo-finite exponentiation}

Fix notation for a non-standard model of $\C$
$${^*\C}=\C^P/\mathcal{D}$$
where $P$ and the ultrafilter $\mathcal{D}$ on $P$ are defined in \ref{1.3}.

Note that by construction $\ZZ\subset {^*\C}$ and this allows us to identify elements $l\in \F_\pp$ which are represented by $l\in \ZZ,$ $0\le l< \pp,$ with $l\in  {^*\C}$ in the theorem below.

\bpk \label{Th1}  {\bf Theorem.} {\em 
 Then there is a place  $\mathcal{I}: \F_\mathfrak{p}\to {^*\bar{\C}}$ .
such that $\mathcal{I}$ maps:

for all $l\in \ZZ$ satisfying $-\lfr \le l\le \lfr$ 
  \be\label{muF}  l\mapsto l\ee

\be\label{ccF}\iota\mapsto
e^{-\frac{\pi i}{4}}\mbox{ and } \ii\mapsto e^{-\frac{\pi i}{2}}  \ee
For all $a\in \Q,$ for all $l\in \ZZ,$ $-\lfr <l\le\lfr,$ 
\be \label{xiqmuF} \epsilon^{\frac{al(\pp-1)}{2\lfr}}\mapsto
e^{-\frac{al\pi i  }{\lfr}}\ee

\be \label{xiiF} \epsilon^{ \frac{al(\pp-1)}{2\lfr\ii}}\mapsto e^{-\frac{al\pi}{\lfr} } .\ee

}
\epk

The proof is by Lemmata below.

\medskip

We consider linear equations of the form 
 $$\sum_{i=1}^k c_i\cdot X_i=1$$
 where  the variables $X_i$ are assumed to be in a specific subset $G$ of the field. A solution $x_1,\ldots,x_k$ is said to be {\bf non-degenerate} if for any proper subset $K\subset \{ 1,\ldots,k\}$
 $$\sum_{i\in K} c_i\cdot x_i\neq 1.$$
 
\bpk\label{L0roots} {\bf Lemma.} {\em There is a function $f: \N\to \N,$  an $\eta\in {^*\N}$ and a highly divisible $\nu\in {^*\N}$   such that
 
 $2\nu \eta|(\pp-1),$ 
 
 for all $n\in \N$ $\nu^n| \eta,$ 
 
 and
 
for all $k\in \N$, for all rational functions $c(X,Y)=\la c_1(X,Y),\ldots,c_{k}(X,Y)\ra,$ for any $0\le l\le\nu$ 

\be\label{Pnu}\left.\begin{array}{lll}
\mbox{any non-degenerate solution }  x_1,\ldots,x_{k}\in \F_\pp \\
\mbox{of }\sum_{i=1}^k c_i(l,\eta)\cdot x_i=1\ \& \ \bigwedge_{i}  x_i^\nu=1\\
\mbox{satisfies }\bigwedge_{i}  x_i^{f(k)}=1.
\end{array}
 \right\rbrace
 \ee

}

{\bf Proof.} Recall that $\F_\pp$ is a field of characteristic 0.
Thus the  well-known Theorem of Mann about linear equations in roots of unity with rational coefficients is applicable. 
A consequence of Mann's Theorem is that there is a function $f$ satisfying (\ref{Pnu}) for any $ \nu, \eta\in \N$ 
  (see \cite{Roots} for this and other consequences).

 We treat the expression $x^\nu$ as an arithmetic function of $x,\nu$ defined in $(\ZZ;+,\cdot, \pp)$ along with the interpretation of the field $\F_\pp$.  

Let $\mathcal{M}_c\subset {^*\N}^2$ be the set of $(\nu,\eta)\in {^*\N}^2$ such that (\ref{Pnu}) holds for  given $k$ and $c(X,Y).$ Clearly
 $\mathcal{M}_c$ is definable in $(\ZZ; +,\cdot,\pp).$
 By the above  consequence of the Mann Theorem $\N^2\subseteq \mathcal{M}_c$ and so
 $$\N^2\subseteq \bigcap_{c}\mathcal{M}_c$$
 where $c$ runs in all $k$-tuples of rational functions $c(X,Y).$ 

Since each $n\in \N$ divides $\pp-1$
it follows that the countable type $$\bigwedge_c (\nu,\eta)\in \mathcal{M}_c\ \& \ 2\nu\eta|(\pp-1)\ \&\ \bigwedge_{n\in \N} n|\nu\ \&\ \nu^n|\eta$$
is consistent, thus has a realisation in the $\aleph_0$-saturated structure $\ZZ.$ $\Box$

\medskip

Below we use notation
$$\ZZ[\lfr]:= \{ l\in  \ZZ: -\lfr\le l\le \lfr\}.$$

Assuming $\lfr <<\pp$ we may equally treat  $\ZZ[\lfr]$ as a subset of $\F_\pp.$

\epk
\bpk \label{cormulambda} {\bf Corollary.} {\em We may assume that  for all $k$ and $c(X,Y)$ (\ref{Pnu}) is satisfied when $\nu:=\lfr$ and, if $\ii^2+1\neq \pp,$ $\eta:=\ii$.  
In particular,} \be \label{ccmax}\begin{array}{lll}
 \lfr^n|\ii,\mbox{ for all }n\in\N\mbox{ and }\\
 \ii^2+1= \pp \   \mbox{ or $\ii$ is transcendental in } \F_\pp \mbox{ over } \ZZ[\lfr]

\end{array}\ee
\epk
\bpk
Set 
$$ 'i\R':=\{ { \frac{\kappa}{\lfr}\vv}:\ -m\lfr/2\le \kappa\le m\lfr/2,\ m\in \N\}\subset \U_{\pp,\lfr}$$
and let
$$'\mathbb{S}':=\exp_\pp('i\R')\subset \F_\pp.$$
Since $\exp_\pp(\vv)=1,$ $'\mathbb{S}'$ is the  group of all the elements $\gamma$ of $\F^\times_\pp$ satisfying $\gamma^{\lfr}=1.$

\epk
\bpk\label{S} {\bf Lemma.} {\em Let $\gamma_1,\ldots, \gamma_n\in {'\mathbb{S}'}$ be multiplicatively independent.
Then $\gamma_1,\ldots, \gamma_n$ are algebraically independent over  $\Q(\ZZ[\lfr], \ii).$ 
}

{\bf Proof.} Suppose not. Then for some $k,$ for some  $c_1,\ldots,c_{k}\in \Q(\ZZ[\lfr], \ii),$ for some monomials
$x_i=x_i(\gamma_1,\ldots, \gamma_n),$ $i=1,\ldots, k,$ the equality $\sum_{i=1}^k c_i\cdot x_i=1$ holds. We  assume that $k$ is minimal with this property. Since $\gamma_1,\ldots, \gamma_n$ are roots of 1 of order dividing  $\lfr$  so are their products $x_1,\ldots,x_{k-1}.$ Since $k$ is minimal, the solution  $x_i(\gamma_1,\ldots, \gamma_n),$ $i=1,\ldots, k,$ of the equation is non-degenerate. Thus the $x_i(\gamma_1,\ldots, \gamma_n)$ satisfy (\ref{Pnu}) and so
are roots of unity of order $\le f(k).$ This would contradict our assumption on multiplicative independence. $\Box$ 
\epk
\bpk\label{3.7} {\bf Corollary.} {\em There is a place
$$\mathcal{I}: \Q(\ZZ[\lfr], \ii, {'\mathbb{S}'})\to {^*\bar{\C}}$$
which satisfies (\ref{muF}),(\ref{ccF}) and (\ref{xiqmuF}). }

{\bf Proof.} Denote $\Q(\ZZ[\lfr], \ii, \sqrt[\infty]{1})$ the extension of
$\Q(\ZZ[\lfr], \ii)$ by roots of $1.$ 

First define the place $\mathcal{I}_0: \Q(\ZZ[\lfr],\sqrt[\infty]{1}))\to {^*\bar{\C}}$
to be defined as the obvious embedding of the subfield of $\F_\pp$ generated by
$\ZZ[\lfr]$ and  $\sqrt[\infty]{1})$ into ${^*\bar{\C}},$ where $\ZZ[\lfr]$ is treated as both subset of $\F_\pp$ and of ${^*\C}.$

In case $\ii$ satisfies the option $\ii^2+1=\pp$ of (\ref{ccmax}), $\ii$ is a square root of $-1$ in $\F_\pp$ and so can be identified with an appropriate element of  $\sqrt[\infty]{1})$ and thus included into the domain of $\mathcal{I}_0.$

In the alternative option $\ii$ is independent over the domain of  
$\mathcal{I}_0$ and so one can extend $\mathcal{I}_0$ to an $\mathcal{I}_1$ sending $\ii$ to $\sqrt{-1}$ of ${^*\C}.$

By Lemma \ref{S} $\mathcal{I}_1$ can be extended to $\mathcal{I}$ as required. $\Box$

\epk

\bpk
Let 
$$'\R':=\{ { \frac{\kappa}{\lfr}\uu}:\ -m\lfr/2\le \kappa\le m\lfr/2,\ m\in \N\}\subset \U_{\pp,\lfr},$$
a subgroup of $\U_{\pp,\lfr},$

and
$$'\R'_+:=\exp_\pp('\R')\subset \F_\pp,$$
a subgroup of $\F^\times_\pp.$

\medskip

 {\bf Remark.} Note that the definition of $\exp_\pp$ in \ref{1.4} depends on the choice of the generator $\epsilon.$ So the above $'\R'_+$ 
depends on $\epsilon$ as well (but $'\mathbb{S}'$ does not).  
\epk
\bpk \label{L1Th1} {\bf Lemma.}
 {\em Let 
$k,l\in \N$  and $\alpha=\alpha(\epsilon)=\epsilon^\frac{\pp-1}{\lfr\ii}.$
Let 
$\bar{M}=
\la M_1,\ldots,M_{k}\ra$ in $\ZZ^{k},$ $-l\lfr <M_i< l\lfr.$ 

Let $\bar{g}(\bar{Z})=\la g_1(Z_1,\ldots,Z_k),\ldots,g_k(Z_1,\ldots, Z_k)\ra$ be a $k$-tuple of rational functions over $\Q(\ZZ[\lfr],\ii)$ and $\bar{c}=\la c_1,\ldots,c_{k}\ra,$ $c_i\in \Q(\Z[\lfr],\ii, {'\mathbb{S}'})$ of the form $$c_i=g_i(\bar{s}),\ \bar{s}:=\la s_1,\ldots, s_k\ra,
\mbox{ for } s_1,\ldots, s_k\in {'\mathbb{S}'}.$$

Consider a non-standard Laurent polynomial in $\F_\pp$ of the form 
 \be\label{P} P_{\bar{c},\bar{M}}(X):=\sum_{i=1}^{k} c_i X^{M_i}-1.\ee
 Then there is a generator $\epsilon$ of $\F_\pp$ such that 
for all $\bar{M}$ and $\bar{s}$ as above 
$$P_{\bar{c},\bar{M}}(\alpha)\neq 0.$$}

{\bf Proof.} For a given $\bar{M}$ and $\bar{c}$,
$P_{\bar{c},\bar{M}}(X)$ has at most $l\lfr$ zeroes. There are at most $(l\lfr)^k$ possible tuples $\bar{c}=\bar{g}(\bar{s})$ and 
$(l\lfr)^{k}$ of tuples $\bar{M},$ so at most $(l\lfr)^{2k}$ polynomials (\ref{P}) altogether, so $(l\lfr)^{2k+1}$ zeroes of the polynomials.
 
On the other hand $\alpha(\epsilon)$ takes any value of primitive root of order $\lfr\ii$ as $\epsilon$ runs through generators of $\F^\times_\pp.$ Thus there are $\varphi(\lfr\ii)$ (the Euler function) such $\alpha$ and the well-known  lower estimate gives us
 $\varphi(\lfr\ii)> \sqrt{\lfr\ii},$ 
 which is bigger than $(l\lfr)^{2k+1}$ by
 (\ref{ccmax}). 
 Thus there is an $\epsilon$ such as $\alpha(\epsilon)$ is as required.  $\Box$ 
\epk

\bpk\label{CorTh1} {\bf Corollary.} {\em There is $\epsilon$ such that $\alpha(\epsilon)$ is not a zero of any $P_{\bar{c},\bar{M}}(X)$ for all $l,$  $k,$  $\bar{M}$ and $\bar{c}\in \Q(\ZZ[\lfr],\ii, {'\mathbb{S}'})^k$ as in \ref{L1Th1}.}

{\bf Proof.} 
Note that the conclusion of \ref{L1Th1} can be restated as a formal statement $\exists \epsilon\, \Phi_{l,\bar{g}}(\alpha(\epsilon), \lfr,\ii,\pp)$ in the language of arithmetic. The statement of (\ref{L1Th1})
readily generalises to the statement that $\alpha$ avoids zeroes of a finite number of polynomials of the form (\ref{P}), by taking the product of the polynomials. 
 Thus the type
$$ \{ \Phi_{l,\bar{g}}(\alpha(\epsilon), \lfr,\ii,\pp):\ l, \bar{g} \mbox{ as in (\ref{L1Th1})} \}$$
in variable $\epsilon$ is consistent.
Since $\ZZ$ is $\aleph_0$-saturated there is an $\epsilon$ such that $\alpha(\epsilon)\in  \F_\pp$
 is not a zero of any polynomial (\ref{P}) $\Box$
\epk
\bpk \label{L2Th1} {\bf Lemma.} {\em There is a generator $\epsilon$ such that  any multiplicatively independent   $\gamma_1,\ldots,\gamma_m\in\, '\R'_+$ are algebraically independent over  $\Q(\ZZ[\lfr],\ii, {'\mathbb{S}'}).$  }

{\bf Proof.} Let $\epsilon$ be as stated in \ref{CorTh1}.
 Suppose $\gamma_1,\ldots,\gamma_m$ are multiplicatively independent and satisfy a polynomial equation $R(\gamma_1,\ldots,\gamma_m)=0.$ The equation can be rewritten as
$$\sum_j c_j \gamma^{\bar{m}_j}=1$$
 for some monomials $\gamma^{\bar{m}_j}=\prod_i \gamma_i^{m_{ji}}$ and   $c_j\in \Q(\lfr,\ii,{'\mathbb{S}'}).$
 
  Clearly, $\gamma^{\bar{m}_j}$ belong to the  group $'\R'_+$ and so 
 $\gamma^{\bar{m}_j}=\alpha(\epsilon)^{M_j},$ for some  $M_j,$ $|M_j|<l\lfr,$
  contradicting our choice of $\epsilon.$ $\Box$ 
\epk
 \bpk\label{3.12} {\bf Corollary.} {\em The place $\mathcal{I}$ of \ref{3.7} can be extended to 
 $$\mathcal{I}: \Q(\ZZ[\lfr], \ii, {'\mathbb{S}'},'\R'_+ )\to {^*\bar{\C}}$$
 which satisfies (\ref{xiiF}).}
 
 {\bf Proof.} 
$\mathcal{I}$ of  (\ref{xiiF}) is an isomorphism of groups. Since the only relations in the language of rings between elelements in $'\R'_+$  are multiplicative relations, $\mathcal{I}$ is a place in the language of rings. $\Box$

\epk

\bpk  Since  ${^*\C}$ is algebraically closed the $\mathcal{I}$ of \ref{3.12} can be extended to a place 
$$\mathcal{I}:\F_\pp\to {^*\bar{\C}}.$$
This finishes the proof of Theorem \ref{Th1}. 

$\Box$
\epk 
\bpk \label{Th3} {\bf Proof of the Main Theorem \ref{MainA}: (\ref{R}) -  (\ref{e1})}.

Define $$\ml_\F:= \mathrm{st}\circ \mathcal{I}.$$
In particular,   taking into account that, for $\mathrm{st}: {^*\C}\to \bar{\C},$
$$\mathrm{st}(e^x)=e^{\mathrm{st}(x)}$$
we get  
$$\ml_\F:\ \ \  \epsilon^{\frac{l}{2\lfr}\uu}\mapsto e^{-\pi\, \mathrm{st}(\frac{l}{\lfr})}; \ \ \  \epsilon^{\frac{l}{2\lfr}\vv}\mapsto e^{-i\pi\, \mathrm{st}(\frac{l}{\lfr})}.$$

 It is easy to check that 
$${'\R'}\cap {'i\R'}=\{ 0\}.$$

For $\frac{l}{\lfr}\in {^*\Q}$ (non-standard rationals), such that $\frac{l}{\lfr}\uu\in {'\R'},$ define  
\be\label{e2I} \ml_\U(\frac{l}{\lfr}\uu):=-2\pi \mathrm{st}(\frac{l}{\lfr}).\ee
This is an additive homomorphism 
$${'\R'}\twoheadrightarrow \R.$$
Then also  $\frac{l}{\lfr}\vv\in {'i\R'},$ and define
\be\label{e1I} \ml_\U(\frac{l}{\lfr}\vv):=-i2\pi \, \mathrm{st}(\frac{l}{\lfr}),\ee
an additive homomorphism 
$${'i\R'}\twoheadrightarrow i\R.$$

This defines $\ml_\U$ on ${'\R'}+ {'i\R'}$ respecting the commutation with $\exp.$ Moreover,
$$\ml_\U: \ {'\R'}+ {'i\R'}\twoheadrightarrow \C.$$

Note that the definition of the homomorphism $\ml_\U$ on $'\R'$ and $'i\R'$ above extends uniquely on their divisible hulls
$$\ml_\U:   \frac{1}{n}{'\R'}\to \R \mbox{ and }
\ml_\U:   \frac{1}{n}{'i\R'}\to i\R,\mbox{ \ \ for }n\in \N$$

The sum of divisible hulls $\mathrm{H}({'\R'})+\mathrm{H}({'i\R'})$ is a divisible subgroups of $\U_{\pp,\lfr}$ and so can be   complemented by a subgroup $\U_{\pp,\lfr}(\infty),$
$$\U_{\pp,\lfr}=\mathrm{H}({'\R'})\dot{+} \mathrm{H}({'i\R'})\dot{+}\U_{\pp,\lfr}(\infty).$$
Define $\ml_\U(u)=\infty$ for all $u\in \U_{\pp,\lfr}(\infty).$ 
Using $\exp_\pp$ and $\exp$ as in  the commuting diagram (\ref{Diag1}) 
this can be  extended to $\ml_\F$ so that 
 (\ref{R}), (\ref{iR}), (\ref{ll}), (\ref{e2}) and (\ref{e1}) are satisfied.
$\Box$ 
\epk
The rest of the Main Theorem will be proved in the next section.
\section{Integration}
Recall that $\lfr=\mu^2$ and $\ii=\iota^2$ elements of ${^*\N}.$
\bpk {\bf Proposition.} {\em Let $a=\frac{d^2}{l^2}$ for some $d,l\in \N.$
\be \label{sumu}\sum_{-\ii\lfr/2a\le n\le \ii\lfr/2a} \exp_\pp(-a\frac{n^2}{2\lfr} \uu)= \mu \frac{\iota\omega}{\sqrt{a}}
\ee

\be \label{sumv}\sum_{-\lfr/2a\le n\le \lfr/2a} \exp_\pp(-a\frac{n^2}{2\lfr} \vv)= \mu \frac{\omega}{\sqrt{a}}
\ee
where $\omega\in \F_\pp$ is a primitive root of $1$ of order $8.$
 and $\mu,\iota$ and $d,l$ should be seen as elements of $\F_\pp$ corresponding to the respective integers.
}

{\bf Proof.}
Let $\nu\in {^*\N}$ be even and $2\nu^2|(\pp-1).$ Let $\xi\in \F_\pp$ be a primitive root of order $2\nu^2,$ equivalently, for some $\epsilon\in \F_\pp$, a primitive root of order $\pp-1,$ 
$$\xi=\epsilon^\frac{\pp-1}{2\nu^2}.$$ 

Writing $n=m\nu+k,$ $0\le m,k<\nu,$ we get

\be\label{sum1} \sum_{0\le n<\nu^2} \xi^{n^2}= \sum_{0\le k<\nu} \xi^{k^2}\cdot\sum_{0\le m<\nu} \xi^{m^2\nu^2+2mk\nu}.\ee
 
Now we use the fact that
$\frac{m^2-m}{2}\in  \ZZ$
and get in $\F_\pp$
$$\sum_{0\le m<\nu} \xi^{m^2\nu^2+2mk\nu}=\sum_{0\le m<\nu} \xi^{m{\nu^2}+2mk\nu}=
\sum_{0\le m<\nu} \xi^{2m\nu(\frac{\nu}{2}+k)}=$$

$$=\lbrace \begin{array}{ll}
{\nu}, \mbox{ if } \frac{\nu}{2}+k\equiv 0\mod \nu,\\
0,\mbox{ otherwise}
\end{array}$$
using that $\xi^{2m\nu(\frac{\nu}{2}+k)}=\xi^{2\nu^2}=1$ in the first line
and that $$\sum_{0\le m<\nu}\zeta^m=0,$$ for $\zeta:=\xi^{2\nu(\frac{\nu}{2}+k)},$
$\zeta^\nu=1$ but 
$\zeta\neq 1,$ in the second line.

Hence in the  sum (\ref{sum1})  only $k=\frac{\nu}{2}$ contributes, and we get

\be\label{xisum} \sum_{0\le n<\nu^2} \xi^{n^2}=
{\nu}\,\xi^\frac{\nu^2}{4}=\nu \epsilon^\frac{\pp-1}{8}=\nu\cdot \omega \ee
for $\omega:=\epsilon^\frac{\pp-1}{8},$ a primitive root of 1 of order $8.$

Note that
$$\sum_{0\le n<\nu^2} \xi^{n^2}=\sum_{-\nu^2/2\le n<\nu^2/2} \xi^{n^2}$$
because of periodicity
$$\xi^{(n+\nu^2)^2}=\xi^{n^2}.$$

Finally, set $$\nu:=\frac{\mu\iota}{\sqrt{a}}\mbox{ for (\ref{sumu}) \ and\ }\ \nu:=\frac{\mu}{\sqrt{a}}\mbox{ for (\ref{sumv})} $$ and (\ref{sumu}) and (\ref{sumv}) follow. $\Box$ 

\epk
\bpk\label{4.2} {\bf Corollary (proof of the Main Theorem, (\ref{integral2}) and (\ref{integral1})) } 
 $$ \ml_\F:\
 \frac{1}{\mu}\sum_{-\ii\lfr/a\le n<\ii\lfr/a} \exp_\pp( a\frac{n^2}{2\lfr}\uu)\ \ \mapsto \ \ \frac{1}{\sqrt{a}}$$
 $$ \ml_\F:\  \frac{1}{\mu}\sum_{-\lfr/a \le n<\lfr/a} \exp_\pp( a\frac{n^2}{2\lfr}\vv)\ \ \mapsto\ \ \frac{ e^\frac{\pi i}{4}}{\sqrt{a}}
$$
This follows from the facts that $\ml_\F(\iota)=e^{-\frac{\pi i}{4}}$ and
that $\ml_\F(\epsilon^\frac{\pp-1}{8})=e^{\frac{\pi i}{4}},$ see (\ref{e2I}) and (\ref{e1I}) together with Theorem \ref{Th1}. 

\epk
\bpk
Note that $\frac{1}{\sqrt{a}}=\int_\R e^{- a x^2} dx,$ the classical Gaussian integral. Analogously, in quantum mechanics $\int_\R e^{ia x^2} dx:= \frac{ e^\frac{\pi i}{4}}{\sqrt{a}},$ although the integral is not classically defined since   $e^{ia x^2}$ is oscillating on the whole of $\R.$ One of the ways of justifying the assignment of the value to the integral expression is by referring to the fact that the respective Fresnel integral $\int_{-A}^A e^{ia x^2} dx$ is well-defined for any $A>0$ and 
$$\lim_{A\to \infty}\int_{-A}^A e^{ia x^2} dx= \frac{ e^\frac{\pi i}{4}}{\sqrt{a}}.$$

\epk
\bpk {\bf The domains of integrations and domains of summation.}

Let, for $l\in \N$
 $$\II_l=\{ n\in \ZZ: \
-l\mu\le n\le l\mu\}\mbox{ and } \II=\bigcup_{l\in \N} \II_l.$$
Let $a=\frac{m}{l}$ and
$$\II_{a,\uu}=\{ n\in \ZZ: 
-\ii\lfr/2a\le n\le \ii\lfr/2a\} \mbox{ and }\II_{a,\vv}=\{ n\in \ZZ: 
-\lfr/2a\le n\le \lfr/2a\}$$
the domains of summations of (\ref{sumu}) and (\ref{sumv}). Clearly, using the assumptions on $\lfr$ and $\mu,$
$$\II\subset  \II_{a,\uu}\mbox{ and }\II\subset \II_{a,\vv}$$
and 
$$\frac{\uu}{\mu}\cdot\II\subset {'\R'}\mbox{ and }\frac{\vv}{\mu}\cdot\II\subset {'i\R'}.$$

The application of $\ml_\U$ defined in \ref{Th3} gives us
$$\ml_\U: \ \frac{\uu}{\mu}\cdot\II_l\twoheadrightarrow \R\cap [-l\pi,l\pi]\mbox{ and } \frac{\uu}{\mu}\cdot\II\twoheadrightarrow \R,$$
$$\ml_\U: \ \frac{\vv}{\mu}\cdot\II_l\twoheadrightarrow i(\R\cap [-l\pi,l\pi])\mbox{ and } \frac{\vv}{\mu}\cdot\II\twoheadrightarrow i\R,$$
that is $\frac{\uu}{\mu}\cdot\II$ can be seen as  a Riemann integration partition of sets $'\R'$ with the infinitesimal mesh $\frac{\uu}{\mu},$
and respectively $\frac{\vv}{\mu}\cdot\II$ in $i\R$ with  mesh $\frac{\vv}{\mu}.$

Note that in (\ref{sumu}) and (\ref{sumv}) 
$$\frac{n^2}{2\lfr}=\frac{1}{2}(\frac{n}{\mu})^2;\ \ \ml_\U(\frac{n^2}{2\lfr}\uu)=-\pi \frac{x^2}{2},\ \ \ml_\U(\frac{n^2}{2\lfr}\vv)=-i\pi \frac{x^2}{2}$$
for $x=\mathrm{st}(\frac{n}{\mu}).$
\epk
\bpk\label{4.5} {\bf Lemma.} {\em For every $l\in \N$
$$\ml_\F:\frac{1}{\mu}\sum_{n\in \II_l} \exp_\pp(-a\frac{n^2}{2\lfr} \uu) \mapsto
\int_{-l}^{l} e^{-ax^2}dx$$
$$\ml_\F:\frac{1}{\mu}\sum_{n\in \II_l} \exp_\pp(-a\frac{n^2}{2\lfr} \vv) \mapsto
\int_{-l}^{l} e^{-iax^2}dx$$
and the integrals are well-defined.}

Proof. Note that by (\ref{xiqmuF}) and (\ref{xiiF}) the map
$\mathcal{I}$ translates the respective elements of the sums from $\F_\pp$ to the elements of the non-standard model ${^*\C}$ of  complex numbers:
$$\mathcal{I}: \ \exp_\pp(-a\frac{n^2}{2\lfr} \uu)\mapsto e^{-a\pi(\frac{n}{\mu})^2}$$
$$\mathcal{I}: \ \exp_\pp(-a\frac{n^2}{2\lfr} \vv)\mapsto e^{-ia\pi(\frac{n}{\mu})^2}$$
and thus $$\mathcal{I}\{\frac{1}{\mu}\sum_{n\in \II_l} \exp_\pp(-a\pi\frac{n^2}{2\lfr} \uu)\}\mbox{ and } \mathcal{I}\{\frac{1}{\mu}\sum_{n\in \II_l} \exp_\pp(-a\pi\frac{n^2}{2\lfr} \vv)\}$$
become non-standard Riemann sums with infinitesimal mesh $\frac{1}{\mu}$
$$\sum_{-l< \frac{n}{\mu}<l} \frac{1}{\mu}e^{-a\pi(\frac{n}{\mu})^2}\mbox{ and } \sum_{-l< \frac{n}{\mu}<l} \frac{1}{\mu}e^{-ia\pi(\frac{n}{\mu})^2}$$
Since the summation is over a compact interval $[-l,l]\subset {^*\R}$ the application of the standard part map gives us (see e.g. the integration via non-standard analysis in \cite{Robinson})
 $$\mathrm{st}(\sum_{-l< \frac{n}{\mu}<l} \frac{1}{\mu}e^{-a\pi(\frac{n}{\mu})^2})=\int_{-l}^l e^{-a\pi x^2}dx\mbox{ and } \mathrm{st}(\sum_{-l< \frac{n}{\mu}<l} \frac{1}{\mu}e^{-ia\pi(\frac{n}{\mu})^2})=\int_{-l}^l e^{-ia\pi x^2}dx$$

\epk
\bpk\label{4.6} {\bf Corollary.} $$\ml_\F:\frac{1}{\mu}\sum_{n\in \II} \exp_\pp(-a\frac{n^2}{2\lfr} \uu) \mapsto
\int_\R e^{-ax^2}dx=\frac{1}{\sqrt{a}}$$
$$\ml_\F:\frac{1}{\mu}\sum_{n\in \II} \exp_\pp(-a\frac{n^2}{2\lfr} \vv) \mapsto
\lim_{l\to \infty}\int_{-l}^{l} e^{-iax^2}dx=\frac{e^\frac{\pi}{4}}{\sqrt{a}}$$

Proof. In both cases the right-hand side is the limit of integrals in \ref{4.5} since $\II=\bigcup_l\II_l.$ In the first case the classical Gaussian  integral over the whole of $\R$ converges and is equal to the limit. In the second case the right-hand side as the limit of the Fresnel integral is well-defined  but the Riemann integral is not. 
\epk
\bpk {\bf Discussion}. The left-hand sides of \ref{4.2} and of \ref{4.6} differ in the domains of summations but the right-had sides are the same. This implies that our definition of $\ml_\F$ is such that
 $$\ml_\F\left(\frac{1}{\mu}\sum_{n\in \II_{a,\uu}\setminus \II} \exp_\pp(-a\frac{n^2}{2\lfr} \uu)\right)=0 \mbox{ and }\ml_\F\left(\frac{1}{\mu}\sum_{n\in \II_{a,\vv}\setminus \II} \exp_\pp(-a\frac{n^2}{2\lfr} \vv)\right)=0$$
In these ``tail domains'' $\II_{a,\uu}\setminus \II$ and $\II_{a,\vv}\setminus \II$ the respective values under exponentiation are very large, non-feasible numbers, and according to the inerpretation of $\F_\pp$ in $\C$ the application of $\exp$ to such values oscillates uncontrollably. R.Feynman intuition was that for this reason the sum should be considered negligible, see e.g. \cite{Feynman}, 2-3.
\epk
 \end{document}